\documentclass[11pt]{amsart}

\theoremstyle{definition}

\theoremstyle{remark}

\begin{document}

\title{A map on the space of rational functions}

\author[Boros]{G. Boros}
\address{Department of Mathematics, Xavier University, New
Orleans, Louisiana 70125}
\email{gboros@xula.edu}

%\author[Briscoe]{S. Briscoe}
%\address{Department of Mathematics,
%Tulane University, New Orleans, LA 70118}
%\email{sage@math.tulane.edu}

\author[Little]{J. Little}
\address{Department of Mathematics and Computer Science,
College of the Holy Cross, Worcester, MA 01610}
\email{little@mathcs.holycross.edu}

%\author[Manna]{D. Manna}
%\address{Department of Mathematics,
%Tulane University, New Orleans, LA 70118}
%\email{manna@math.tulane.edu}

\author[Moll]{V. Moll}
\address{Department of Mathematics,
Tulane University, New Orleans, LA 70118}
\email{vhm@math.tulane.edu}

\author[Mosteig]{E. Mosteig}
\address{Department of Mathematics,
Loyola Marymount University, Los Angeles, CA 90045}
\email{emosteig@lmu.edu}

\author[Stanley]{R. Stanley}
\address{Department of Mathematics,
Massachusetts Institute of Technology, Cambridge, MA 02139-4307}
\email{rstan@math.mit.edu}

%    General info
\subjclass{Primary 33}

\date{\today}


\begin{abstract}
We describe dynamical properties of a map ${\mathfrak{F}}$ 
defined on the space of rational functions.
The fixed points of ${\mathfrak{F}}$ are  classified and the 
long time behavior of a subclass is described in terms 
of Eulerian polynomials. 
\end{abstract}

\maketitle

%\textwidth=6.5in
%\textheight=9.0in
%\topmargin=-\headheight
%\oddsidemargin=0in
%\evensidemargin=0in
%\def\Tilde{\char126\relax}

\newcommand{\nn}{\nonumber}
\newcommand{\ba}{\begin{eqnarray}}
\newcommand{\ea}{\end{eqnarray}}
\newcommand{\A}{{\mathfrak{A}}}
\newcommand{\Am}{{\mathfrak{A}_{m}}}
\newcommand{\B}{{\mathfrak{B}}}
\newcommand{\C}{{\mathbb{C}}}
\newcommand{\F}{{\mathfrak{F}}}
\newcommand{\R}{{\mathfrak{R}}}
\newcommand{\N}{{\mathbb{N}}}
\newcommand{\Cn}{{\mathfrak{C}_{n}}}
\newcommand{\ift}{\int_{0}^{\infty}}
\newcommand{\no}{\noindent}

\newtheorem{Definition}{\bf Definition}[section]
\newtheorem{Thm}[Definition]{\bf Theorem}
\newtheorem{Lem}[Definition]{\bf Lemma}
\newtheorem{Cor}[Definition]{\bf Corollary}
\newtheorem{Note}[Definition]{\bf Note}
\newtheorem{Prop}[Definition]{\bf Proposition}
\numberwithin{equation}{section}

\section{Introduction} \label{intro}
\setcounter{equation}{0}

Let $\R$ denote the space of rational functions with complex coefficients. The
Taylor expansion at $x =0$ of $R \in \R$ is written as
\ba
R(x) & = & \sum_{n \gg -\infty} a_{n}x^{n}
\ea
\no
where $n \gg -\infty$ denotes the fact that the coefficients vanish
for large negative $n$.  We
consider the map $\F: \R \to \R$ defined by
\ba
\F(R(x)) & = & \sum_{n \gg -\infty} a_{2n+1}x^{n}.
\ea
\no
The map $\F$ can be given explicitly by
\ba
\F(R(x)) & = & \frac{R(\sqrt{x}) - R( - \sqrt{x}) }{2 \sqrt{x}}
\ea
\no
and it appeared in this form in the description of a general procedure for the
exact integration of rational functions \cite{bomolan2}. The splitting of
an arbitrary function $R$ into its even and odd parts $R(x) = R_{e}(x) +
R_{o}(x)$ yields
\ba
\ift R(x) \, dx & = & \ift R_{e}(x) \,dx + \ift R_{o}(x) \, dx.
\ea
\no
The integral of the even part can be analyzed with the methods described in
\cite{bomolan2}, and the integral of the odd part can be transformed by
$ x \mapsto \sqrt{x}$ to produce
\ba
\ift R(x) \, dx & = & \ift R_{e}(x) \, dx + \tfrac{1}{2}
\ift \F(R(x)) \, dx.
\ea

\medskip

Here we consider dynamical properties of the map $\F: \R
\rightarrow \R$. Section \ref{kernel} characterizes rational
functions $R$ for which the orbit \ba \text{Orb}(R) & := & \left\{
\F^{(k)}(R(x)): \, k \in \N \right\} \ea \no ends at the fixed
point $0 \in \R$. Section \ref{special} describes
the dynamics of a special class of functions with all their poles
restricted to the unit circle. We establish explicit formulas for
the asymptotic behavior of their orbits, expressed in terms of Eulerian
polynomials $A_{m}(x)$ defined by the generating function 
\ba 
\frac{1-x}{1-x \,
\text{exp}[ \lambda (1-x) ]} & = & \sum_{m=0}^{\infty} A_{m}(x)
\frac{\lambda^{m}}{m!}. \label{eulerdef} \ea 
\no 
The proof only employs the classical result \cite{comtet}, page 243:
\ba
A_{m}(x) & = & (1-x)^{m+1} \sum_{k=0}^{\infty} k^{m}x^{k}.
\label{eulerdef0}
\ea
\no
Section \ref{fixedpoints} contains a description of
all the fixed points of $\F$.

\medskip

The map $\F$ can be supplemented by
\ba
\mathfrak{E}(R(x)) & := & \frac{R( \sqrt{x}) + R(- \sqrt{x})}{2}
\ea
\no
for which similar results can be established. See \cite{simu2002a}
for details. These transformations
can be written as
$$
\mathfrak{F}(R(x)) = \frac{R(x^{1/2}) +
 \omega_{2}R( \omega_{2}x^{1/2})}{2 x^{1/2}} \text{ and }
\mathfrak{E}(R(x)) = \frac{R(x^{1/2}) +
 R(\omega_{2}x^{1/2})}{2}
$$
\no
where $\omega_{2} = -1$. The extensions to higher degree will be
considered in future work. For instance, in the case of degree $3$
we consider the maps:
\ba
{\mathfrak{T}}_{1}(R(x)) & := & \frac{1}{3} \left[ R( x^{1/3} ) +
R( \omega_{3}x^{1/3} ) +
R( \omega_{3}^{2}x^{1/3} ) \right] \nn \\
{\mathfrak{T}}_{2}(R(x)) & := & \frac{1}{3x^{1/3}} \left[ R( x^{1/3} ) +
\omega_{3}^{2}R( \omega_{3}x^{1/3} ) +
\omega_{3}R( \omega_{3}^{2}x^{1/3} ) \right] \nn \\
{\mathfrak{T}}_{3}(R(x)) & := & \frac{1}{3x^{2/3}} \left[ R( x^{1/3} ) +
\omega_{3}R( \omega_{3}x^{1/3} ) +
\omega_{3}^{2}R( \omega_{3}^{2}x^{1/3} ) \right] \nn
\ea
\no
where $\omega_{3} = e^{2 \pi i /3}$. These maps correspond to subseries 
of the expansion  of $R$ taken along a fixed class modulo $3$. \\

\no
{\bf Notation}: 
\no
For $n \in \mathbb{N}, \, m_{1}, \cdots, m_{n}$ are odd integers.  \\

\no
$L = (m_{1}m_{2} \cdots m_{n})^{-1} $ and 
$m^{*}  = m_{1} + m_{2} + \cdots + m_{n} -2.$  \\

\no
For $j \in \mathbb{Z}$ and $n \in \mathbb{N}$
\ba
Q_{n}(x) & = & \prod_{k=1}^{n} (x^{m_{k}}-1) \label{qn}
\ea
\no
and
\ba
R_{j,n}(x) & =  & \frac{x^{j}}{Q_{n}(x)}. \label{rjn}
\ea

\medskip

\section{The kernel of the iterates} \label{kernel}
\setcounter{equation}{0}

Even rational functions can be characterized as the elements of the kernel
of $\F$. In this section we characterize those functions that vanish
precisely after $n$ applications of $\F$.  \\

The sets
\ba
K_{n} & := & \text{Ker } {{\mathfrak{F}}}^{(n)}
      \;  =
\; \left\{ S \in \mathfrak{R}: \; {{\mathfrak{F}}}^{(n)}(S) = 0 \; \right\}
\nn
\ea
\no
form a nested sequence of vector spaces. We now describe the class of 
functions $J_{n} := K_{n} \, - \, K_{n-1}$, i.e.
those functions that vanish after precisely $n$ applications of the
map $\mathfrak{F}$. In particular we show that $J_{n}$ is not empty.

\medskip

The decomposition of $S \in \mathfrak{R}$ into its even and odd
parts can be expressed as \ba S(x) & = & S_{1,1}(x^{2}) + x
S_{2,1}(x^{2}) \ea \no where $S_{1,1}, \, S_{2,1} \in
\mathfrak{R}$. This decomposition applied to $S_{1,1}, \, S_{2,1}$
yields 
\ba S(x) & = & S_{1,2}(x^{4}) + x S_{2,2}(x^{4}) + x^{2}
S_{3,2}(x^{4}) + x^{3} S_{4,2}(x^{4}). \nn 
\ea
\no
Iterating this procedure produces a general decomposition.  \\

\begin{Lem}
Given $n \in \mathbb{N}$ and $S \in \mathfrak{R}$ there is a unique
set of rational functions $\{ S_{j,n}: \; 0 \leq j \leq 2^n-1 \}$ such that
\ba
S(x) & = &  \sum_{j=1}^{2^{n}} x^{j-1} S_{j,n}(x^{2^{n}}).
\label{ratdecomp}
\ea
\end{Lem}
\begin{proof}
Split the sum in
\ba
S(x) & = & \sum_{k \gg - \infty} a_{k}x^{k} \label{expanrat}
\ea
\no
according to the residue of $k$ modulo
$2^{n}$.
\end{proof}

\medskip

We now show that the functions $S \in J_{n}$ are precisely those for which
$S_{2^{n},n} = 0$ and $S_{2^{n-1},n} \neq 0$.  This generalizes the
case $n=1$ that states that $\F(S) = 0$ precisely when
$S$ is even.

\medskip

\begin{Thm}
\label{thm-preper}
The rational functions that vanish after precisely $n$ applications of $\F$
are those of the form
\ba
S(x) & = &  \sum_{j=1}^{2^{n}-1} x^{j-1} S_{j,n}(x^{2^{n}}),
\ea
\no
where $S_{1,n}, \cdots, S_{2^{n}-2,n}$ are arbitrary rational functions and
$S_{2^{n}-1,n} \neq 0$.
\end{Thm}
\begin{proof}
A direct calculation from (\ref{ratdecomp}) shows that for $1 \leq k$, 
\ba
{\F}^{(k)}(S)(x) & = & \sum_{j=1}^{2^{n-k}} x^{j-1} S_{2^{k}j,n}(x^{2^{n-k}}).
\ea
\no
The statement now follows from
\ba
{\F}^{(n-1)}(S)(x) & = & S_{2^{n-1},n}(x^{2}) + x S_{2^{n},n}(x^{2})
\ea
\no
and
\ba
{\F}^{(n)}(S)(x) & = & S_{2^{n},n}(x).
\ea
\end{proof}

\no
{\bf Note}. We now state the result of Theorem \ref{thm-preper} in the
language of dynamical systems. Let $f: X \to X$ be a map  on a set $X$ and
$x_{0} \in X$ be a fixed point of $f$, i.e. $f(x_{0}) = x_{0}$. We say that a
sequence
$\{ x_{1}, \, x_{2}, \, \cdots, \, x_{n} \}$ of elements of $X$ is
a {\em prefixed sequence of length} $n$ {\em attached to} $x_{0}$
if $ x_{i+1} = f(x_{i})$ and
$f(x_{n}) = x_{0}$. Theorem \ref{thm-preper} states that $\F$ admits prefixed
sequences attached to $0$ of arbitrary length.

\medskip

\section{The dynamics of a special class} \label{special}
\setcounter{equation}{0}

The asymptotic behavior of $\F$ can be described in complete detail for
rational functions in the class
\begin{equation}
\Cn := \Cn \left(m_{1}, \ldots, m_{n} \right) = \left\{
\frac{P(x)}{Q_{n}(x)}: \, P \text{ is a Laurent polynomial}  \; \right\},
\end{equation}
\no
where $m_{1}, \cdots, m_{n}$ are {\em odd positive integers} and
$Q_{n}(x)$ is defined in (\ref{qn}). 
A {\em Laurent polynomial} is a
rational function of the form
$ a_{-k} x^{-k} + a_{-k+1}x^{-k+1} + \cdots + a_{j-1}x^{j-1} + a_{j}x^{j}$
with $k, \, j \in \mathbb{N}$. \\

\medskip

The case $n=1$ has been described in \cite{bjm}. The results are expressed
in terms of the function
\ba
\gamma_{m}(j) & = & m \left\lfloor \frac{j}{2} \right\rfloor -
\frac{1}{2}(m-1)(j-1) \label{gamma} \\
 & = & \begin{cases}
        (m-1+j)/2  & \text{ if } j \text{ is even} \\
        (j-1)/2   & \text{ if } j \text{ is odd.}
        \end{cases}
\nn
\ea

\medskip

\begin{Thm}
Let $m$ be an odd positive integer. Then  \\

\no
1) For $j \in \mathbb{N}$
\ba
\F \left( \frac  {x^{j}}   {x^{m}-1} \right) &  = &
\frac{   x^{\gamma_{m}(j)}    }{x^{m}-1}.
\ea
\no
Thus, the study of the dynamics of $\F$ on ${\mathfrak{C}}_{1}$ is reduced to 
that of $\gamma_{m}$ on $\mathbb{Z}$.

\medskip

\no
2) The iterates $\left\{ \gamma_{m}^{(p)}(j): \, p=0, \, 1, \, 2, \, \cdots
\right\}$ reach the set
\ba
\Am & := & \left\{ 0, \, 1, \, 2, \, \cdots, m-2 \right\}
\ea
\no
or the fixed points $-1$ and $m-1$ in a finite number of 
steps. Moreover, $\Am$ is invariant under the action
of $\gamma_{m}$.  This action partitions $\Am$ into orbits.

\medskip

\no
3) The inverse of the restriction of $\gamma_{m}$ to $\Am$ is given by
\ba
\delta_{m}(k) & = & \begin{cases}
                    2k+1  & \text{    if }  0 \leq k \leq (m-2)/2  \\
                    2k+1-m  & \text{ if }  (m-1)/2 \leq k \leq m-2.
                     \end{cases}
\ea
\end{Thm}

\no
\begin{Note}
The explicit form of $\delta_{m}$ permits the explicit computation
of the orbits of $\gamma_{m}$ on the invariant set $\Am$. For example, if
$m$ is prime then every orbit of $\gamma_{m}$ is of length $\text{Ord}(2;m)$.
In particular, there is a single orbit if and only if $2$ is a primitive
root modulo $m$. \\
\end{Note}

The iterates $\gamma_{m}^{(p)}(j)$ can be characterized by the
congruence (\ref{congr1}) below. This
will be used in the determination of the limiting behavior of the
iterates of $\F$ below.  The proof of this congruence and the numerical 
and symbolic evidence of the asymptotic behavior of the iterates 
of $\F$ on $R_{j,n}(x)$
were part of a SIMU 2002 project. Details will appear in \cite{simu2002a}. \\

\begin{Lem}
\label{lemmacongr}
Let $m$ be an odd positive number, $0 \leq j < m$, and $p \in \mathbb{N}$.
The unique solution of
\ba
2^{p}(x + 1) & \equiv & j+1 \quad \text{ mod } m
\label{congr1}
\ea
\no
in $0 \leq x < m$ is given by $x = \gamma_{m}^{(p)}(j)$.

\end{Lem}
\begin{proof}
The proof is by induction on $p$. Note first that the solution of the
congruence (\ref{congr1}) is unique because $\text{gcd}(m,2^{p})=1$.
The base case ($p=0$) is
$$x + 1 \equiv j+1 \quad \text{ mod }m$$
\no
with unique solution $j = \gamma_{m}^{(0)}(j)$.  To complete the
inductive step observe that
\ba
2^{p+1} \left( \gamma_{m}^{(p+1)}(j) + 1 \right)  & = &
2^{p+1} \left( m \left\lfloor \frac{\gamma_{m}^{(p)}(j)}{2} \right\rfloor
- \frac{1}{2}(m-1) ( \gamma_{m}^{(p)}(j) -1 ) +  1 \right) \nn \\
& = & 2^{p+1}m \left\lfloor \frac{\gamma_{m}(j)}{2} \right\rfloor -
2^{p}m \left( \gamma_{m}^{(p)}(j) - 1 \right) +
  2^{p} \left( \gamma_{m}^{(p)}(j) +1 \right) \nn \\
 & \equiv & j+1 \quad \text{ mod }m.  \nn
\ea
\end{proof}

\medskip

The map $\F$ has a very rich dynamical structure, even in the case $n=1$. \\

\no
\begin{Thm}
Let $r \in \mathbb{N}$. Then $\F$ has at least one periodic orbit of
length $r$.
\end{Thm}
\begin{proof}
For $r \in \mathbb{N}$ define $m = 2^{r}-1$. Then the orbit of $1/(x^{m}-1)$
under $\F$ is
$$
\frac{1}{x^{m}-1} \mapsto  \frac{x^{2^{r-1}-1}}{x^{m}-1}  \mapsto
\frac{x^{2^{r-2}-1}}{x^{m}-1} \mapsto \, \cdots  \,
\mapsto \frac{x}{x^{m}-1} \mapsto
\frac{1}{x^{m}-1}.$$
\end{proof}

\medskip

We now consider the properties of the class $\Cn$ for $n >  1$. \\

\begin{Lem}
\label{lemmapres}
The class $\Cn$ is invariant under the action of $\F$.
\end{Lem}
\begin{proof}
We show that $\F(R_{n,j}(x)) \in \Cn$.
The linearity of $\F$ yields the result.  \\

Introduce the notation
\ba
S_{l} := \sum x^{(m_{i_{1}} + m_{i_{2}} + \cdots + m_{i_{l}}) /2}
\ea
\no
where the sum is taken over all subsets of $\{ m_{1}, \cdots, m_{n} \}$
containing $l$ elements, the empty sum giving $S_{0}=1$.

Now
\ba
\prod_{k=1}^{n} \left( x^{m_{k}/2} + 1 \right) & = & 1 + \sum_{l=1}^{n}S_{l}
\ea
\no
and
\ba
\prod_{k=1}^{n} \left( x^{m_{k}/2} - 1 \right) & = & (-1)^{n}
+ \sum_{l=1}^{n}(-1)^{n-l} S_{l},
\ea
\no
so that
\ba
\F ( R_{j,n}(x))  & = & \frac{x^{(j-1)/2}}{2Q_{n}(x)}
\left\{ \sum_{l=0}^{n} S_{l} + (-1)^{j+1} \sum_{l=0}^{n} (-1)^{l}S_{l}
\right\} \nn \\
 & =: & \frac{x^{(j-1)/2}}{2 Q_{n}(x)} \times S. \nn
\ea
\no
The sum $S$ contains all the terms $S_{l}$ with $l$ of the opposite parity
of $j$ so that $x^{(j-1)/2} \times S$ is a polynomial.
\end{proof}

\medskip

We now consider the orbit of a general rational function $R$ in
the class $\Cn$. The case $n=2$ illustrates the general situation. A direct
calculation shows that \\

\ba
\F \left( \frac{x^{j}}{(x^{m_{1}}-1)(x^{m_{2}}-1)} \right) & = &
\frac{x^{(j-1+m_{1})/2} + x^{(j-1+m_{2})/2}}{(x^{m_{1}}-1)(x^{m_{2}}-1)}
\text{ if } j \text{ is even}  \nn
\ea

\medskip

\no
and

\ba
\F \left( \frac{x^{j}}{(x^{m_{1}}-1)(x^{m_{2}}-1)} \right) & = &
\frac{x^{(j-1)/2} + x^{(j-1+m_{1}+m_{2})/2}}{(x^{m_{1}}-1)(x^{m_{2}}-1)}
\text{ if } j \text{ is odd}.  \nn
\ea

\no
Thus $\F$ preserves the denominator of $R$ (as was shown
in Lemma \ref{lemmapres}) and each monomial of $R$ yields two monomials. We
now show that the exponents of these monomials are always bounded.  \\

\begin{Prop}
\label{bounds} Let $m^{*} = \sum\limits_{k=1}^{n} m_{k} -2$ and define
\ba \A & := & \left\{ \frac{P}{Q_{n}} \in \Cn: \; P \text{ is a
polynomial with } \text{deg}(P) \leq m^{*} \right\}. \nn \ea \no
Then $\mathfrak{A}$ is invariant under $\F$, and every orbit
starting at $R \in \Cn$ reaches it in a finite number of steps.
\end{Prop}
\begin{proof}
Let $R = P/Q_{n} \in \Cn$ and write $\F(R) = P_{1}/Q_{n}$. The largest
exponent in $P_{1}$, call it $j$, appears from the sum $S_{n}$.  The inequality
$j > m^{*}$ implies $\tfrac{1}{2}(j-1 + m_{1} + \cdots + m_{n} ) < j$. The 
case $j < 0$ is similar.
\end{proof}

\medskip

Now we consider the asymptotic behavior of the iterates of $\F$ starting
at $R \in \Cn$. 
This is expressed in terms of the Eulerian polynomials $A_{m}(x)$
defined in (\ref{eulerdef}). The discussion is divided into two cases 
depending whether the number 
\ba 
d & := &
{\rm{gcd}}(m_{1}, \, m_{2}, \cdots, m_{n}) \label{gcd1} \ea \no
is $1$ or not.  In 
Theorem \ref{gcdone} we prove that if $d=1$ then 
\ba
\F^{(p)}(R(x)) & \sim & \frac{2^{(n-1)p}}{(n-1)! \, m_{1} \cdots m_{n}}  
\, \frac{A_{n-1}(x)}{(1-x)^{n}}.
\ea
\no
The case $d > 1$ is described in Theorem \ref{gcdbig1}. We prove that
the sequence of iterates of $\F$ applied to the function
$x^{j}/Q_{n}(x)$ taken along a fixed residue class
(defined in (\ref{rhodef})) has an
asymptotic behavior as in the case $d=1$, with limit points
expressed in terms of Eulerian polynomials. \\
% as \ba \F^{(p \rho(j) +
%q)} \left( \frac{x^{j}}{Q_{n}(x)} \right) & \sim & 2^{(n-1)(p
%\rho(j) +q)} \frac{x^{\gamma}}{(x^{d}-1)^{n}} \sum_{l=0}^{n-1}
%\beta_{l} (x^{d}-1)^{n-1-l} A_{l}(x^{d}) \nn \ea \no
%where $\gamma$ and $\beta_{l}$ are constants.  \\

This suggests the existence of
an arbitrary number of limit points, but we have not ruled out the possibility
that all these could coincide.   \\

\medskip

The proof employs the observation that if
\begin{equation}
R(x) = \sum_{k} f(k) x^{k}
\end{equation}
\no
is the expansion of $R$, then
\begin{equation}
\F^{(p)}(R)(x) = \sum_{k} f(2^{p}(k+1)-1) x^{k}.
\end{equation}

\medskip

\begin{Thm}
\label{gcdone}
Let $R \in \A$ and suppose $d=1$. Then
\ba
\lim\limits_{p \to \infty} \frac{\F^{(p)}(R)(x)}{2^{(n-1)p}} & = &
\frac{L \, A_{n-1}(x)}{(1-x)^{n} \, (n-1)!}
\ea
\no
where $L = 1/(m_{1} \cdots m_{n})$.
\end{Thm}
\begin{proof}
Since $d=1$, the rational function
\ba
R(x) & = & \frac{P(x)}{(x^{m_{1}}-1) \cdots (x^{m_{n}}-1)}
\ea
\no
has a pole of order $n$ at $x=1$ and poles of order less than $n$ at the
other zeros of $Q_{n}(x)$, all of which are roots of unity. Thus the partial
fraction expansion of $R$ has the form
\ba
R(x) & = & \frac{L}{(1-x)^{n}} + G(x)
\ea
\no
where the term $G(x)$ contains all the terms of order lower than $n$. Hence
\ba
R(x) & = & L \sum_{k \geq 0} \binom{n+k-1}{n-1} x^{k} + G(x)
\ea
\no
where the coefficient of $x^{k}$ in $G(x)$ is $O(k^{n-2})$.  Thus
\ba
\F^{(p)}(R)(x) & = & L \sum_{k} \binom{n+ 2^{p}k+2^{p}-2}{n-1}x^{k} + 
\F^{(p)}(G)(x) \nn  \\
& = & \frac{L \, 2^{(n-1)p}}{(n-1)!} \sum_{k} k^{n-1}x^{k} + O(2^{(n-2)p})
\nn
\ea
\no
where we have used
\ba
\binom{n+2^{p}k + 2^{p}-2}{n-1}  & = & 
             \frac{2^{(n-1)p}}{(n-1)!} k^{n-1} + O(2^{(n-2)p}) \nn
\ea
\no
as $p \to \infty$. The result now follows from 
(\ref{eulerdef0}).
\end{proof}

The analysis of the case $d > 1$ involves
the function 
\ba 
\label{rhodef}
\rho(i) & = & \text{ Ord} \left( 2;
\frac{d}{\text{gcd}(i+1,d)} \right) 
\ea 
\no 
with $d$ as in
(\ref{gcd1}). The function $\rho$ appears in the
dynamics of $\gamma_{m}$: for $0 \leq j \leq m-2$ the length of the 
orbit containing $j$ is $\rho(j)$. See \cite{bjm} for details. \\

Introduce the notation 
\ba
E_{i} & := & \{ j \in \Bbb{Z}: \text{ there exists } r \in \Bbb{N}: 
\gamma_{d}^{(r)}(j) = i \; \} 
\ea
\no
for the backward orbit of $\gamma_{d}$ and let 
\ba
A & := & E_{-1} \cup E_{d-1} 
\ea
\no
be the integers that eventually are mapped to the fixed points of 
$\gamma_{d}$.  \\

\begin{Thm}
\label{gcdbig1}

Let $n \in \Bbb{N}$.  \\

\no
a) If $j \in A$ then 
\ba
\lim\limits_{p \to \infty} \frac{\F^{(p)} ( R_{j,n}(x) ) }{2^{(n-1)p}} 
& = & 
\frac{L d^{n} A_{n-1}(x^{d})}{(n-1)! \, (x^{d}-1)^{n} }. 
\nn
\ea

\medskip

\no
b) If $j \not \in A$ choose $r^{*}$ such that 
$$ j^{*} := \gamma_{d}^{(r^{*})}(j)   \in \{ 0, \, 1, \, \cdots , \, d-1 \}.$$ 
Then for any $q \in \{ 0, \, 1, \ldots, \rho(j^{*}) -1 \}$ 
$$\lim\limits_{p \to \infty}
\frac{
\F^{(p \rho(j^{*}) + q)}(R_{j,n}(x))}{
2^{(n-1)(p \rho(j^{*}) + q)}} 
$$
$$
= 
\frac{L x^{\gamma}}{(n-1)! (x^{d}-1)^{n}} 
\times \sum\limits_{l=0}\limits^{n-1}
(-1)^{l} \binom{n-1}{l} \left( \gamma +1 \right)^{l} d^{n-l}
(x^{d}-1)^{l} A_{n-1-l}(x^{d}) $$
\no
where $\gamma = \gamma_{d}^{(q)}(j^{*})$.
\end{Thm}
\begin{proof}
Since $\text{gcd}(m_{1},m_{2}, \cdots, m_{n}) = d$, the function
$R_{j,n}(x)$ has a pole of order $n$ at the $d$-th roots of unity
and a pole of strictly lower order at all other zeros of
$Q_{n}(x)$. Thus its partial fraction decomposition has
the  form
\ba
R_{j,n}(x) & = & 
\frac{C_{0,j}}{(x-1)^{n}} + \frac{C_{1,j}}{(x-\omega_{1})^{n}} +
\cdots + \frac{C_{d-1,j}}{(x - \omega_{d-1})^{n}} +
\text{ lower order terms} \nn
\ea
\no
where $\omega_{l} = \text{exp}(2 \pi i l/d), \; 1 \leq l \leq d-1$.
The coefficients $C_{l,j}$ are given by
\ba
C_{l,j}  =  \lim\limits_{x \to \omega_{l}} \frac{(x- \omega_{l})^{n} x^{j}}
{Q_{n}(x)}
  =  L \, \omega_{l}^{n+j}. \nn
\ea
\no
Thus
\ba
R_{j,n}(x) & = & \sum_{l=0}^{d-1} \frac{C_{l,j}}{(x - \omega_{j})^{n}}
+ G(x), \label{extraterm}
\ea
\no
where the error term $G(x)$ includes the polar 
parts of all poles of order less than
$n$. Then
\ba
R_{j,n}(x)  & = & \sum_{l=0}^{d-1} C_{l,j} \frac{(-1)^{n}}{\omega_{l}^{n}}
\sum_{k=0}^{\infty} \binom{n+k-1}{n-1} \frac{x^{k}}{\omega_{l}^{k}} + G(x) \nn \\
 & =& (-1)^{n}L \sum_{k=0}^{\infty} \sum_{l=0}^{d-1} \binom{n+k-1}{n-1}
\omega_{l}^{j-k} x^{k} + G(x). \nn
\ea
\no
Therefore \\

$
\lefteqn{
\F^{(r)} (R_{j,n}(x)) = 
}
$
$$
(-1)^{n}L \sum_{k=0}^{\infty}
\left( \sum_{l=0}^{d-1} \omega_{l}^{j-2^{r}(k+1)+1} \right)
\binom{n+2^{r}(k+1)-2}{n-1} x^{k} + \F^{(r)}(G(x)).
$$
\no
Now
\ba
\sum_{l=0}^{d-1} \omega_{l}^{q} & = & \begin{cases}
                                  0  \quad \text{ if }   d \text{ does not divide }q
                                   \\
                                  d  \quad \text{ if }   d \text{ divides } q,   \\
                                      \end{cases}
\nn \ea \no so the only values of $k$ that contribute to the sum
are those for which 
\ba
j+1 & \equiv  & 2^{r}(k+1) \quad \text{
mod  }d. \nn
\ea
\no
The discussion is divided into two cases depending on whether some iterate
of $j$ reaches one of the fixed points or not. \\

\no
{\bf Case 1}: Suppose $j \in A$ and let $r^{*} \in \Bbb{N}$ be such that 
$\gamma_{d}^{(r^{*})}(j) = -1$. Then 
Lemma \ref{lemmacongr} shows that  
$k = \gamma_{d}^{(r+r^{*})}(j) + N d = Nd -1$.  Thus  \\

$
\lefteqn{
{\F}^{(r+r^{*})}( R_{j,n}(x))  = 
}
$
$$
 (-1)^{n} Ld 
\sum_{N=1}^{\infty} \binom{n-2 + 2^{r+r^{*}} Nd }{n-1} x^{Nd-1} + 
{\F}^{(r+r^{*})}(G(x)). 
$$
\no
In the limit as $r \to \infty$ the binomial coefficient is asymptotic to 
$$ \frac{2^{(r+r^{*})(n-1)} N^{n-1} d^{n-1}}{(n-1)!}. $$
\no
The number $r^{*}$ is fixed, so we have 
\ba
\lim\limits_{r \to \infty} 
\frac{{\F}^{(r+r^{*})}(R_{j,n}(x)) }{2^{(r+r^{*})(n-1)}} & = & 
\frac{(-1)^{n} L d^{n} A_{n-1}(x^{d})}{(n-1)! \, (1- x^{d})^{n}}
\ea
\no
as stated.   \\

A similar argument shows that the same result is true if $j$ lies on the 
backward orbit of the second fixed point. \\

\no
{\bf Case 2}: Now assume $j \not \in A$. Lemma \ref{lemmacongr} shows that 
$$k = \gamma_{d}^{(r+r^{*})}(j) + N d =  \gamma_{d}^{(r)} 
\left( \gamma_{d}^{(r^{*})}(j) \right) + Nd. $$ 
\no

\no
Thus  \\

$
\lefteqn{
{{\F}^{(r+r^{*})} ( R_{j,n}(x) ) } =   
}
$
$$
(-1)^{n}Ld \sum_{N=0}^{\infty} \binom{n-2 + 2^{r+r^{*}} ( Nd +  
\gamma_{d}^{(r+r^{*})}(j) + 1 )}{n-1} x^{Nd + \gamma_{d}^{(r+r^{*})}(j) }. 
$$ 

\no
In the limit as $r \to \infty$, the binomial coefficient is asymptotic to 
$$
\frac{2^{(r+r^{*})(n-1)} (Nd + \gamma_{1} + 1)^{n-1}}{(n-1)!}
$$
\no
where $\gamma_{1} = \gamma^{(r+r^{*})}(j)$. Thus 
$$
\frac{ {\F}^{(r+r^{*})}(R_{j,n}(x))}{2^{(r+r^{*})(n-1)}}  =  
\frac{(-1)^{n} Ld x^{\gamma_{1}}}{(n-1)!} 
\sum_{N=0}^{\infty} (Nd + \gamma_{1} + 1)^{n-1} \, x^{Nd}  + o(1). 
$$
\no
as $r \to \infty$. Now
\ba
\sum_{N=0}^{\infty} (Nd + \gamma_{1} + 1 )^{n-1} \, x^{Nd}  & = & 
\sum_{N=0}^{\infty} x^{Nd} \, 
\sum_{l=0}^{n-1} \binom{n-1}{l} ( 1+ \gamma_{1})^{n-1-l} N^{l} d^{l}  \nn \\
& = & \sum_{l=0}^{n-1} \binom{n-1}{l} ( 1+ \gamma_{1})^{n-1-l}  d^{l} 
\frac{A_{l}(x^{d})}{(1-x^{d})^{l+1}}.  \nn
\ea
\no
As a function of $r$,
$$ \gamma_{1}  =  \gamma^{(r+r^{*})}(j) = \gamma^{(r)}(j^{*}) $$
\no
has period $\rho(j^{*})$. Write $r = p \rho(j^{*}) + q$ with 
$0 \leq q \leq \rho(j^{*})-1$ and replace
$\gamma_{d}^{(r)}(j^{*})$ by $\gamma_{d}^{(q)}(j^{*})$ to obtain  \\
$$
\frac{ {\F}^{(p \rho(j^{*}) + q+r^{*})}(R_{j,n}(x))}
{2^{(p \rho(j^{*}) + q+r^{*})(n-1)}}  =  
$$
$$
\frac{(-1)^{n} Ld x^{\gamma_{d}^{(q)}(j^{*})}}{(n-1)!} 
\sum_{l=0}^{n-1} \binom{n-1}{l} ( 1 + \gamma_{d}^{(q)}(j^{*}) )^{l}  \,
\times \frac{d^{n-1-l} A_{n-1-l}(x^{d})}{(1-x^{d})^{n-l}} + o(1).
$$
\no
To conclude the proof, observe that as $q$ runs over the set  of residues 
modulo $\rho(j^{*})$, so does $q + r^{*}$.
\end{proof}

\section{The fixed points of $\F$} \label{fixedpoints}
\setcounter{equation}{0}

A formal power series argument shows that any rational function
$R$ fixed by $\F$ must have an expansion of the form 
\ba 
x R(x) & = & 
c + \sum_{n=0}^{\infty} a_{n} \varphi(x^{2n+1}),
\label{exprat} 
\ea 
\no 
where 
\ba \label{expphi}\varphi(x) & = &
\sum_{n=0}^{\infty} x^{2^{k}} = x + x^{2} + x^{4} + x^{8} +
\cdots. \ea 
\no
In particular, $R$ has at most a simple pole at the
origin and if \ba R(x) & = & \sum_{n \geq -1} f(n)x^{n} \nn \ea
\no is such a fixed point, then $f(2n) = f(n)$ for $n \geq 0$.
Thus the problem of finding fixed points of $\F$ is reduced to
finding sequences $\{ a_{n} \}$ for
which (\ref{exprat}) is a rational function. \\

The class of functions discussed in Section \ref{special} yields examples of
fixed points. Let $m$ be an odd positive integer. Then $m-1$ is 
fixed by $\gamma_{m}$, so
$R_{1,m-1}(x) = x^{m-1}/(x^{m}-1)$ is fixed by $\F$.  This example
can be obtained by a different approach. First observe that if
$\F(R) = R$ and $r$ is any odd positive integer, then the function
\ba
\B_{r}(R(x))  & = & x^{r-1}R(x^{r})
\ea
\no
is also fixed by $\F$. The function $g(x) = 1/(x-1)$ is fixed by $\F$, so
that $R_{1,m-1}(x) = \B_{m}(g(x)) $ is also fixed. \\

The description of all the fixed points of $\F$ requires the
notion of {\em cyclotomic cosets}:  given $n, \, r \in \mathbb{N}$
with $r$ odd and  $ 0 \leq n \leq r-1$, the set 
\ba 
C_{r,n} & = &
\{ 2^{s}n \, \text{ mod } \, r: \; s \in \mathbb{Z} \} 
\ea 
\no 
is the $2$-cyclotomic coset of $n$ mod $r$. Observe that $C_{r,n}$ is 
a finite set. With $\lambda$
a fixed primitive $r$-th root of unity, define \ba f_{r,n}(x) & =
& \sum_{m \in C_{r,n}} \frac{\lambda^{m}}{1 - \lambda^{m} x}. \ea

The partial fraction decomposition of the fixed point $\B_{m}(g(x))$
can be decomposed into a sum of rational functions each fixed by $\F$.
For example, consider $\B_{7}(g(x)) = x^{6}/(x^{7}-1)$ and let
$\lambda = \text{exp}(2 \pi i/7)$ be a primitive $7$-th root of
unity. Then
\ba
7 \B_{7}(g(x)) & = & \sum_{k=0}^{6} \frac{\lambda^{k}}{1 - \lambda^{k}x}  \\
          & = & \frac{1}{1-x} +
    \left( \frac{\lambda}{1 - \lambda x} + \frac{\lambda^{2}}{1- \lambda^{2}x}
         + \frac{\lambda^{4}}{1 - \lambda^{4}x} \right) + \nn \\
    &   & \left( \frac{\lambda^{3}}{1 - \lambda^{3} x} + \frac{\lambda^{5}}{1-
    \lambda^{5}x}
         + \frac{\lambda^{6}}{1 - \lambda^{6}x} \right), \nn
\ea \no where each of the sets of terms grouped together is a
rational function fixed by $\F$. In the notation introduced above,
this decomposition is \ba 7 \B_{7}(g(x)) & = & f_{7,0}(x) +
f_{7,1}(x) + f_{7,3}(x). \ea

\medskip

We now classify all the fixed points of $\F$. \\

\medskip

\begin{Thm}
\label{fixedpt}
A rational function is fixed by $\F$ if and only if it is a linear combination of
$\tfrac{1}{x}$ and
the functions $f_{r,n}(x)$ for $r$ odd and $0 \leq n \leq r-1$.
\end{Thm}
\begin{proof}
The identity
\ba
\F \left( \frac{\lambda}{1 - \lambda x} \right) & = & \frac{\lambda^{2}}{1 -
\lambda^{2}x} \nn
\ea
\no
shows that $\F$ fixes the $f_{r,n}$ because the squaring map $S: z \mapsto z^{2}$
permutes the
values $\lambda^{m}$ for $m \in C_{r,n}$.  \\

We first establish the converse under the assumption that the
poles of $R$ are simple. The final step of the proof consists of
checking
that this condition holds for any fixed point of $\F$.  \\

Let $R$ be a rational function, with simple poles,  that is fixed by $\F$. The
partial fraction decomposition of $R$ is \ba R(x) & = &
\frac{c}{x} + \sum_{j=1}^{J} \frac{\alpha_{j} \lambda_{j}} {1 -
\lambda_{j}x}, \label{unique} \ea \no which is unique up to order.
Apply $\F$ to produce \ba R(x) & = & \frac{c}{x} + \sum_{j=1}^{J}
\frac{\alpha_{j} \lambda_{j}^{2}} {1 - \lambda_{j}^{2}x}. \nn \ea
\no The uniqueness of (\ref{unique}) shows that \ba \{
\lambda_{1}, \cdots , \lambda_{J} \} & = & \{ \lambda_{1}^{2},
\cdots , \lambda_{J}^{2} \}, \ea \no that is, the set $\Lambda =
\{ \lambda_{1}, \cdots, \lambda_{J} \}$ is permuted by the
squaring map. We conclude that the set $\{ \lambda_{k}^{2^{n}}: \,
n \in \mathbb{N} \}$ is a finite set (a subset of $\Lambda$) and
so every $\lambda_{k}$ is an $r_{k}$-th
root of unity, for some odd positive integer $r_{k}$. \\

Now group terms in the sum (\ref{unique}) according to the orbits
of the squaring map $S$ on the set $\Lambda$. The coefficient
$\alpha_{j}$ in (\ref{unique}) must be constant along each orbit, and 
moreover, the orbit of $\lambda_{k}$ under $S$ is precisely the
set $\{ \lambda_k^{m}: \, m \in C_{n,r_{k}} \}$ for some $n$.
Therefore $R(x)$ can be decomposed as a  linear combination
of the required form. \\

The next result concludes the proof of the theorem.

\begin{Prop}
\label{simple} Let $R(x)$ be a rational function that can be
expressed in the form \ba x R(x) & = & c + \sum_{n=0}^{\infty}
a_{n} \varphi(x^{2n+1}) \ea \no where $\varphi$ is given in
(\ref{expphi}).
 Then
the poles of $R$ must be simple.
\end{Prop}

\medskip

The following result, used in the proof of Proposition \ref{simple}, is
demonstrated in \cite{stan}, page 202. \\

\begin{Lem}
Let $q_1, q_2, \dots, q_d$ be a fixed sequence of complex numbers,
$d \ge 1$, and $q_d \not= 0$. The following conditions on a
function $f: \N \to \C$ are equivalent:
\begin{enumerate}
\item $\sum_{n \ge 0} f(n)x^n = \frac{P(x)}{Q(x)}$
where, $Q(x) = 1 + q_1 x + q_2 x^2 + q_3 x^3 + \cdots + q_d x^d$.
\item For  $n \gg 0$,
$$f(n) = \sum_{i=1}^k P_i(n) \lambda_i^n,$$
where $1+ q_1 x + q_2 x^2 + q_3 x^3 + \cdots + q_d x^d=
\prod_{i=1}^k (1 - \lambda_i x)^{d_i}$, the $\lambda_i$'s are
distinct, and $P_i(n)$ is a polynomial in $n$ of degree less than
$d_i$.
\end{enumerate}
\end{Lem}

\medskip

\no
{\em Proof of Prop.} \ref{simple}: Assume
$P(x),Q(x)$ are relatively prime, and that $f(n)$ is the
generating function for $P(x)/Q(x)$ written in the form promised
by the lemma above (for $n \gg 0$).  Since $f(n) = f(2n)$,
$$Q(x) = \prod_{i=1}^k (1 - \lambda_i x)^{d_i}
= \prod_{i=1}^k (1- \lambda_i^2x)^{e_i},$$ so
$$\{ \lambda_1, \dots, \lambda_k \} = \{ \lambda_1^2, \dots, \lambda_k^2
\}.$$  As in the proof of Theorem \ref{fixedpt}, we conclude that
each $\lambda_i$ is a primitive $r_i$-th root of unity for some
positive integers $r_1, \dots, r_k$.

Let $M = $ lcm$(r_1,
\dots, r_k)$ and for $a\in \N$, define 
$$R_a = \{ m \in \N : m
\equiv a \; \text{ mod } M \}$$
\no
and $f_a = f \Big{|}_{R_{a}}$ to be
the restriction of the
function $f : \N \to \C$ to the set $R_a$. Then
$$f_a(a+jM) = \sum_{i=1}^k P_i(a+jM) \lambda_i^{a+jM} =
\sum_{i=1}^k P_i(a+jM) \lambda_i^{a},$$ so each $f_a$ has a
representation as a polynomial in the variable $j$ since
$\lambda_i^a$ is constant on the set $R_a$. We denote the natural
extension of this map to an element of $\C[j]$ by $F_a$. Note that
the restriction of $F_a$ to $\N$ will not be $f$ in general. Our
goal is to prove that each $F_a$ is a constant function, with
corresponding constant denoted by $c_a$.
 Once this is shown, we have
 $$\frac{P(x)}{Q(x)} = \sum_{a=1}^M c_a \sum_{j=0}^{\infty}
 x^{a+jM} = \sum_{a=1}^M \frac{c_ax^a}{1-x^M}
 ,$$
so $P(x)/Q(x)$ is a rational function with only simple poles,
as desired. \\

It remains to show  that each polynomial map $F_a: \C \to \C$ is a
constant function.   For each positive integer $n$, define $S_n =
\{ 2^t n : t \in \N  \}$.  We say that $a$ has an {\em infinite
cross-section} if $R_a \cap S_n$ is an infinite set for some $n \in
\N$. We proceed by considering two cases, depending on whether $a$
has an infinite cross-section or not. \\

\no
{\bf Case 1}:  Suppose $a$ has an infinite cross-section, i.e., $R_a
\cap S_n$ is an infinite set.  Since $f(s) = f(2s)$ for all $s\in
\N$ , $F_a$ is constant on $R_a \cap S_n$.  Since $R_a \cap S_n$
is an infinite set, $F_a$ is a constant polynomial. \\

\no
{\bf Case 2}:  Suppose $a$ does not have an infinite cross-section,
i.e.,  $R_a \cap S_n$ is finite for all positive integers $n$.
Then $R_a \cap S_n$ must be nonempty for infinitely many values of
$n$.  Since there are only finitely many distinct sets of the form
$R_b$, it follows that for  each $S_n$, there exists $b \in \N$
such that $R_b \cap S_n$ is infinite.  Moreover, since there are
only finitely many choices for $R_b$, there is at least one $b \in
\N$ such  that there exist infinitely many values of $n$
 where $R_a \cap S_n$ is nonempty and $R_b \cap S_n$ is
infinite.  Since $b$ has an infinite cross-section, an application
of Case 1 demonstrates that the restriction of $f$ to $R_b$ is the
constant function $c_b$. Since $f$ is constant on each $S_n$, the
restriction of $f$ to $S_n$ is the constant $c_b$. Thus
$F_a$ achieves the value $c_b$ infinitely many times, and so $F_a$
must be a constant polynomial.
\end{proof}

\no
The proof of Theorem \ref{fixedpt} is complete.

\medskip

\no
\begin{Note}
For each fixed point $R_{f}$ we construct the rational function
\ba
R(x) & = &  \sum_{j=1}^{2^{n}-1} x^{j-1} R_{j,n}(x^{2^{n}}) - R_{f}(x),
\ea
\no
where $R_{1,n}, \cdots, R_{2^{n}-2,n}$ are arbitrary rational functions and
$R_{2^{n}-1,n} \neq 0$. Theorem \ref{thm-preper} shows that $R$ is the
general form of a prefixed sequence of length $n$ attached to
$R_{f}$.

\end{Note}

\medskip

\no
{\bf Acknowledgments}. 
The third author acknowledges the partial support of 
NSF-DMS 0070567, the last author acknowledges the partial support of 
NSF-DMS 9988459.


\medskip
\begin{thebibliography}{99}

\bibitem{bomolan2}
Boros, G. - Moll, V.: Landen transformations and the integration of
rational functions.
Math. Comp. {\bf 71}, 2002, 649-668.

\bibitem{bjm}
Boros, G. - Joyce, M. - Moll, V.: {\em A transformation on the space of
rational functions}. Elemente der Mathematik {\bf 57}, 2002, 1-11.

\bibitem{simu2002}
Briscoe, S. - Jimenez, L. - Medina, L.: {\em Asymptotics of  a
transformation on the space of rational functions}. SIMU 2002,
Report.

\bibitem{simu2002a}
Briscoe, S. - Jimenez, L. - Manna, D. - Medina, L. - 
Moll, V.: {\em The dynamics of a
transformation on the space of rational functions}. In preparation. 

\bibitem{comtet}
Comtet, L.: {\em Advanced Combinatorics}. Revised and enlarged edition.
D. Reidel Publ. Co., Boston, 1974.

\bibitem{mollnot}
Moll, V.: {\em The evaluation of integrals: a personal story}. Notices AMS,
March 2002, {\bf 49}, 311-317.

\bibitem{stan1}
Stanley, R.: {\em Enumerative Combinatorics. Volume} 1. Cambridge
Studies in Advanced Mathematics {\bf 49}. Cambridge University
Press, 1997.


\bibitem{stan}
Stanley, R.: {\em Enumerative Combinatorics. Volume} 2. Cambridge
Studies in Advanced Mathematics {\bf 62}. Cambridge University
Press, 1999.


\end{thebibliography}
\end{document}